\documentclass[12pt]{article}

\newbox\squ  % box character for ends of proofs
\setbox\squ=\hbox{\vrule width.3pt
            \vbox{\hrule height.3pt width.4em\kern1ex\hrule height.3pt}%
            \vrule width.3pt}
\def\endproof{%
  \ifmmode\eqno\copy\squ\smallskip\else{\unskip\nobreak\hfil%
    \penalty50\hskip2em\hbox{}\nobreak\hfil\copy\squ
    \parfillskip=0pt \finalhyphendemerits=0\penalty-100\smallskip}
  \fi}
\usepackage{amsmath}
\usepackage{amssymb}

\newcommand{\non}{\nonumber}
\newcommand{\wt}{\widetilde}
\newcommand{\wh}{\widehat}
\newcommand{\ot}{\otimes}
\newcommand{\g}{\mathfrak g}
\newcommand{\ts}{\,}
\newcommand{\U}{ {\rm U}}
\newcommand{\CC}{ {\rm C}}
\newcommand{\Z}{ {\rm Z}}
\newcommand{\R}{ {\rm R}}
\newcommand{\J}{ {\rm J}}
\newcommand{\Y}{ {\rm Y}}
\newcommand{\C}{\mathbb{C}}
\newcommand{\ZZ}{\mathbb{Z}}

\newcommand{\sgn}{ {\rm sgn}\ts}
\newcommand{\Norm}{ {\rm Norm}\ts}
\newcommand{\h}{\mathfrak h}
\newcommand{\gl}{\mathfrak{gl}}
\newcommand{\oa}{\mathfrak{o}}

\newcommand{\Proof}{\noindent{\it Proof.}\ \ }

\newtheorem{thm}{Theorem}[section]
\newtheorem{prop}[thm]{Proposition}
\newtheorem{lem}[thm]{Lemma}
\newtheorem{cor}[thm]{Corollary}

\newcommand{\bth}{\begin{thm}}
\renewcommand{\eth}{\end{thm}}
\newcommand{\bpr}{\begin{prop}}
\newcommand{\epr}{\end{prop}}
\newcommand{\ble}{\begin{lem}}
\newcommand{\ele}{\end{lem}}
\newcommand{\bco}{\begin{cor}}
\newcommand{\eco}{\end{cor}}
\newcommand{\beq}{\begin{equation}}
\newcommand{\eeq}{\end{equation}}

\begin{document}

\title{\Large\bf Weight bases of Gelfand--Tsetlin type for\\
representations of classical Lie algebras}
\author{{\sc A. I. Molev}\\[15pt]
School of Mathematics and Statistics\\
University of Sydney,
NSW 2006, Australia\\
alexm@maths.usyd.edu.au}%\\[30pt]
%Research Report 99--??}

\date{} % August 1999
\maketitle

\vspace{7 mm}

\begin{abstract}
This paper completes a series devoted to
explicit constructions of finite-dimensional irreducible
representations of the classical Lie algebras.
Here the case of odd orthogonal	Lie algebras (of type $B$)
is considered (two previous papers dealt with $C$ and $D$ types).
A weight basis for each representation
of the Lie algebra $\oa(2n+1)$ is constructed.
The basis vectors are parametrized by
Gelfand--Tsetlin-type patterns. Explicit formulas for the matrix elements 
of generators of $\oa(2n+1)$ in this basis are given. 
The construction is based on
the representation theory of the Yangians.
A similar approach is applied to the $A$ type case where the well-known
formulas due to Gelfand and Tsetlin are reproduced.

\end{abstract}

\newpage

\section{Introduction}
\setcounter{equation}{0}

In this paper we  give an explicit construction
of each finite-dimensional irreducible representation $V$
of an odd orthogonal Lie algebra $\oa(2n+1)$ (i.e. a simple complex
Lie algebra of type $B$). A weight basis in $V$ is
obtained by the application of certain elements of the enveloping algebra
(the lowering operators) to the highest weight vector.
Explicit formulas for the matrix elements 
of generators of the Lie algebra $\oa(2n+1)$ in this
basis are given. 
We follow an approach
applied in the previous papers \cite{m:br} and \cite{m:wb}
where similar results were obtained for the $C$ and $D$ type
Lie algebras. 
We also reproduce a slightly modified
version of the well-known construction of the
Gelfand--Tsetlin bases for the $A$ type Lie algebras.

Let $\g_n$ denote the rank $n$ simple Lie algebra of type
$A,B,C,$ or $D$. The restriction of a finite-dimensional
irreducible representation $V$ of $\g_n$ to the subalgebra
$\g_{n-1}$ is multiplicity-free for the $A$ type case,
and it is not necessarily so for the $B,C,D$ types.
Gelfand and Tsetlin \cite{gt:fdu} used the chain
of subalgebras
\beq\label{chainInt}
\g_1\subset\g_2\subset\cdots\subset\g_n
\end{equation}
to parametrize basis vectors in $V$ and give formulas 
for the matrix elements
of generators for the $A$ type case. Different approaches
to derive these formulas are used e.g. in
\cite{z:cg, nm:ol, z:gz, g:me, nt:yg, m:gt}.

Analogous results for representations of
the orthogonal Lie algebra are obtained by
Gelfand and Tsetlin in~\cite{gt:fdo};
see also \cite{ph:lr, w:ro, g:wc}. 
Here the chain \eqref{chainInt} is replaced
with the one which involves both the odd and even
orthogonal Lie algebras. However, the corresponding
basis vectors lose the weight property, i.e. they are
not eigenvectors for the elements of the Cartan subalgebra.
To get a weight basis we propose to
use the chain \eqref{chainInt} for the $B,C,D$ types as well.
We ``separate" the multiplicities occurring in the reduction 
$\g_n\downarrow\g_{n-1}$ by applying the representation theory 
of the Yangians. Namely, the subspace $V^+_{\mu}$
of $\g_{n-1}$-highest vectors of weight $\mu$ in $V$ possesses
a natural structure of a representation of the
{\it twisted Yangian\/} $\Y^+(2)$ or $\Y^-(2)$, in the orthogonal
and symplectic case, respectively. The twisted Yangians are
introduced and studied by Olshanski~\cite{o:ty}; see also \cite{mno:yc}
for a detailed exposition. The action of $\Y^{\pm}(2)$
in the space $V^+_{\mu}$ arises from his {\it centralizer
construction\/} \cite{o:ty}.
Finite-dimensional irreducible representations
of the  twisted Yangians
are classified in \cite{m:fd}. In particular, it turns
out that the representation $V^+_{\mu}$ of $\Y^{\pm}(2)$
can be extended to a larger algebra, the {\it Yangian\/} $\Y(2)$
for the Lie algebra $\gl(2)$. 
The algebra $\Y(2)$ and its representations are very well studied;
see \cite{t:im}, \cite{cp:yr}. In particular,
a large class of representation of $\Y(2)$ admits
Gelfand--Tsetlin-type bases associated with the
inclusion $\Y(1)\subset\Y(2)$; see \cite{m:gt, nt:ry}.
This allows us to get a natural basis in the space $V^+_{\mu}$,
and then by induction to get a basis in the entire space $V$.

Note that in the case of $C$ or $D$ type the
$\Y(2)$-module $V^+_{\mu}$ is irreducible while
in the $B$ type case it is a direct sum
of two irreducible submodules. This does not lead, however,
to major differences in the constructions, and the final
formulas are similar in all the three cases.

Our calculations of the matrix elements of the generators of $\g_n$
are based on the relationship between
the twisted Yangian
$\Y^{\pm}(2)$ and the {\it transvector
algebra\/} $\Z(\g_n,\g_{n-1})$
(it is also called the {\it Mickelsson algebra\/} or 
$S$-{\it algebra\/}). It is generated by the {\it raising\/}
and {\it lowering operators\/} which preserve the subspace $V^+$
of $\g_{n-1}$-highest vectors in $V$.
The algebraic structure of the
transvector algebras is studied in detail in \cite{z:it}
with the use of the extremal projections for
reductive Lie algebras \cite{ast:po}.
We construct an algebra homomorphism $\Y^{\pm}(2)\to \Z(\g_n,\g_{n-1})$
which allows us to express the generators of the twisted Yangian,
as operators in $V^+_{\mu}$,
in terms of the raising and lowering operators. This plays
a key role in the calculation of the matrix elements
of the generators of $\g_n$ in the basis of $V$.

Explicit combinatorial constructions of the
fundamental representations of the symplectic
and odd orthogonal Lie algebras were recently given by
Donnelly~\cite{d:ec}. He also showed that in the symplectic case
the basis of  \cite{m:br} for the fundamental representations
coincides, up to a scaling, with a basis of his \cite{d:ecf}.
It is likely that a similar connection exists in the
odd orthogonal case.

\section{Gelfand--Tsetlin basis for $\gl(n)$}\label{sec:gtb}
\setcounter{equation}{0}

Let $E_{ij}$, $i,j=1,\dots,n$ denote the standard basis of the general
linear Lie algebra
$\g_n=\gl(n)$ over the field of complex numbers.
The subalgebra $\g_{n-1}$ is spanned by
the basis elements $E_{ij}$ with $i,j=1,\dots,n-1$. Denote by $\h=\h_n$ the diagonal
Cartan subalgebra in $\g_n$. The elements $E_{11}, \dots,E_{nn}$ form
a basis of $\h$.

Finite-dimensional irreducible representations of $\g_n$
are in a one-to-one correspondence with $n$-tuples
of complex numbers $\lambda=(\lambda_1,\dots,\lambda_n)$ 
such that
\beq%\label{cond}
\lambda_i-\lambda_{i+1}\in\ZZ_+\qquad\text{for}\quad i=1,\dots,n-1. 
\non
\end{equation}
Such an $n$-tuple $\lambda$ is called the highest weight
of the corresponding representation which
we shall denote by $L(\lambda)$.
It contains a unique, up to a multiple, nonzero vector $\xi$
(the highest vector) such that
$
E_{ii}\ts\xi=\lambda_i\ts\xi
$
for $i=1,\dots,n$ and
$
E_{ij}\ts\xi=0
$
for $1\leq i<j\leq n$.
Denote by $L(\lambda)^+$ the subspace of $\g_{n-1}$-highest vectors
in $L(\lambda)$:
\begin{equation}%\label{asubhv}
L(\lambda)^+=\{\eta\in L(\lambda)\ |\ E_{ij}\ts \eta=0,
\qquad 1\leq i<j<n\}.
\non
\end{equation}
Given a $\g_{n-1}$-highest weight
$\mu=(\mu_1,\dots,\mu_{n-1})$ we denote by $L(\lambda)^+_{\mu}$
the corresponding weight subspace in $L(\lambda)^+$:
\begin{equation}\label{ahvwmu}
L(\lambda)^+_{\mu}=\{\eta\in L(\lambda)^+\ |\ E_{ii}\ts\eta=
\mu_i\ts\eta,\qquad i=1,\dots,n-1\}.
\non
\end{equation}
It is well-known \cite{z:cg} that 
the space
$L(\lambda)^+_{\mu}$
is either trivial or one-dimensional. 
Moreover, $\dim L(\lambda)^+_{\mu}=1$
if and only if
\beq\label{amulam}
\lambda_i-\mu_i\in\ZZ_+ 	\quad\text{and}\quad \mu_i-\lambda_{i+1}\in
\ZZ_+ 	\qquad\text{for}\quad i=1,\dots,n-1.
\end{equation}
In other words, the restriction of $L(\lambda)$ to the subalgebra
$\g_{n-1}$ is multiplicity-free:
\beq\label{abranch}
L(\lambda)|^{}_{\g_{n-1}}\simeq \bigoplus L'(\mu),
\non
\end{equation}
where $L'(\mu)$ is the irreducible $\g_{n-1}$-module with the highest
weight $\mu$ satisfying the conditions \eqref{amulam}.
A parameterization of basis vectors in $L(\lambda)$ is obtained by
using its further restrictions to the subalgebras of the chain
\eqref{chainInt}.
%$$
%\g_1\subset\g_2\subset\cdots\subset\g_{n-1}\subset\g_n.
%$$
A {\it Gelfand--Tsetlin pattern\/} $\Lambda$ associated with
$\lambda$ is an array of row vectors
\begin{align}
&\qquad\lambda^{}_{n1}\qquad\lambda^{}_{n2}
\qquad\qquad\cdots\qquad\qquad\lambda^{}_{nn}\non\\
&\qquad\qquad\lambda^{}_{n-1,1}\qquad\ \ \cdots\ \ 
\ \ \qquad\lambda^{}_{n-1,n-1}\non\\
&\quad\qquad\qquad\cdots\qquad\cdots\qquad\cdots\non\\
&\quad\qquad\qquad\qquad\lambda^{}_{21}\qquad\lambda^{}_{22}\non\\
&\quad\qquad\qquad\qquad\qquad\lambda^{}_{11}  \non
\end{align}
such that the upper row coincides with $\lambda$ and 
the following conditions hold
\beq\label{aconl}
\lambda^{}_{ki}-\lambda^{}_{k-1,i}\in\ZZ_+,\qquad  
\lambda^{}_{k-1,i}-\lambda^{}_{k,i+1}\in\ZZ_+,\qquad 
%\qquad\text{for}\quad
i=1,\dots,k-1
\end{equation}
for each $k=2,\dots,n$. 

\medskip
\noindent
{\it Remark.} If the highest weight $\lambda$ is a partition then
there is a natural
bijection between the patterns associated with $\lambda$ and
semistandard $\lambda$-tableaux with entries in $\{1,\dots,n\}$.
A pattern can be viewed as a sequence of partitions
\beq\label{aseq}
\lambda^{(1)}\subseteq \lambda^{(2)}\subseteq\cdots 
\subseteq \lambda^{(n)}=\lambda, \non
\end{equation}
with 
$\lambda^{(k)}=(\lambda^{}_{k1},\dots,\lambda^{}_{kk})$. 
Conditions \eqref{aconl} mean that the skew diagram
$\lambda^{(k)}/\lambda^{(k-1)}$ is a horizontal strip;
see e.g. \cite{ma:sf}. \endproof

Let us set
%\beq\label{alla}
$l^{}_{ki}=\lambda^{}_{ki}-i+1$.
%\end{equation}

\bth\label{thm:abasis} There exists a basis 
$\{\xi^{}_{\Lambda}\}$ in $L(\lambda)$ parametrized by all
patterns $\Lambda$ such that the action
of generators of $\g_n$ is given by the formulas
\begin{align}\label{aekk}
E_{kk}\ts \xi^{}_{\Lambda}&=\left(\sum_{i=1}^k\lambda^{}_{ki}
-\sum_{i=1}^{k-1}\lambda^{}_{k-1,i}\right)
\xi^{}_{\Lambda},
%\non
\\
\label{aekk+1}
E_{k,k+1}\ts \xi^{}_{\Lambda}
&=-\sum_{i=1}^k \frac{(l^{}_{ki}-l^{}_{k+1,1})\cdots (l^{}_{ki}-l^{}_{k+1,k+1})}
{(l^{}_{ki}-l^{}_{k1})\cdots \wedge\cdots(l^{}_{ki}-l^{}_{kk})}
\ts
\xi^{}_{\Lambda+\delta^{}_{ki}},\\
\label{aek+1k}
E_{k+1,k}\ts \xi^{}_{\Lambda}
&=\sum_{i=1}^k \frac{(l^{}_{ki}-l^{}_{k-1,1})\cdots (l^{}_{ki}-l^{}_{k-1,k-1})}
{(l^{}_{ki}-l^{}_{k1})\cdots \wedge\cdots(l^{}_{ki}-l^{}_{kk})}
\ts
\xi^{}_{\Lambda-\delta^{}_{ki}}.
\end{align}
The arrays $\Lambda\pm\delta^{}_{ki}$
are obtained from $\Lambda$ by replacing $\lambda^{}_{ki}$
by $\lambda^{}_{ki}\pm1$. It is supposed
that $\xi^{}_{\Lambda}=0$ if the array $\Lambda$ is not a pattern;
the symbol $\wedge$ indicates that the zero factor in the denominator
is skipped.
\eth

\Proof Consider
the extension 
of the universal enveloping algebra $\U(\g_n)$
\begin{equation}\label{aext}
\U'(\g_n)=\U(\g_n)\ot_{\U(\h)} \R(\h),
\non
\end{equation}
where $\R(\h)$ is the field of fractions of the commutative algebra
$\U(\h)$. Let $\J$
denote the left ideal in $\U'(\g_n)$ generated by 
the elements $E_{ij}$ with $1\leq i<j<n$.
Introduce the normalizer of $\J$ in $\U'(\g_n)$:
\beq\label{anorm}
\Norm \J=\{x\in \U'(\g_n)\ |\ \J\ts x\subseteq \J\}.
\non
\end{equation}
Note that $\J$ is a two-sided ideal in the algebra $\Norm \J$. We define
the {\it transvector algebra\/} $\Z(\g_n,\g_{n-1})$ as the
quotient
\begin{equation}\label{atra}
\Z(\g_n,\g_{n-1})=\Norm \J/\J;
\non
\end{equation}
see \cite{z:gz,z:it}. Equivalently, $\Z(\g_n,\g_{n-1})$ can be defined
by using the {\it extremal projection\/} $p=p_{n-1}$ for the Lie algebra
$\g_{n-1}$ \cite{ast:po}. 
The projection $p$ is, up to a factor from $\R(\h_{n-1})$,
a unique element
of an extension of $\U'(\g_{n-1})$ to an algebra of formal series,
satisfying the conditions
\beq\label{aep}
E_{ij}\ts p=p\ts E_{ji}=0\qquad\text{for}\quad 1\leq i<j\leq n-1.
\end{equation}
The element
$p$ is of zero weight (with respect to the adjoint action of $\h_{n-1}$)
and it can be normalized to satisfy the condition
$
p^2=p.
$
The projection $p$ is a well-defined operator in the
quotient $\U'(\g_n)/\J$ which allows one to naturally
identify the transvector algebra $\Z(\g_n,\g_{n-1})$ 
with the image of $p$ \cite{z:it}:
\beq\label{azp}
\Z(\g_n,\g_{n-1})=p\left(\U'(\g_n)/\J\right).
\non
\end{equation}
An analog of the Poincar\'e--Birkhoff--Witt theorem holds for
the algebra $\Z(\g_n,\g_{n-1})$ so that ordered
monomials in the elements
$pE_{in}$ and $pE_{ni}$ with $i=1,\dots,n-1  $
form a basis of
$\Z(\g_n,\g_{n-1})$ as a left or right $\R(\h)$-module \cite{z:it}.
These elements are called the {\it raising\/} and {\it lowering\/} 
{\it operators\/} and can be given by
\begin{align}%\label{apin}
pE_{in}&=\sum_{i>i_1>\cdots>i_s\geq 1}
E_{ii_1}E_{i_1i_2}\cdots E_{i_{s-1}i_s}E_{i_sn}
\frac{1}{(h_i-h_{i_1})\cdots (h_i-h_{i_s})},
\non\\
\label{apni}
pE_{ni}&=\sum_{i<i_1<\cdots<i_s<n}
E_{i_1i}E_{i_2i_1}\cdots E_{i_si_{s-1}}E_{ni_s}
\frac{1}{(h_i-h_{i_1})\cdots (h_i-h_{i_s})},
\end{align}
where $s=0,1,\dots$ and $h_i=E_{ii}-i+1$. 
We shall also use normalized operators defined by 
\begin{align}%%\label{alow}
z_{in}&=pE_{in}\ts
(h_i-h_{i-1})\cdots (h_i-h_{1}),
\non\\
\label{alow2}
z_{ni}&=pE_{ni}\ts
(h_i-h_{i+1})\cdots (h_i-h_{n-1}).
\end{align}
These can be viewed as elements of the enveloping algebra
$\U(\g_n)$; cf. \cite{z:cg, nm:ol}.
We have the following relations \cite{z:gz, z:it}
\begin{align}\label{acom}
z_{ni}z_{nj}&=z_{nj}z_{ni}\qquad\text{for all}\quad i,j,\\
\label{acom2}
z_{in}z_{nj}&=z_{nj}z_{in}\qquad\text{for}\quad i\ne j.
\end{align}
Indeed, assume that $i<j$. Then \eqref{aep} and
\eqref{apni} imply that in $\Z(\g_n,\g_{n-1})$
\beq%\label{acomp} 
pE_{ni}\ts pE_{nj}=pE_{ni} E_{nj}, \qquad 
pE_{nj}\ts pE_{ni}=pE_{ni} E_{nj}\ts \frac{h_i-h_j+1}{h_i-h_j}.
\non
\end{equation}
Now \eqref{acom} follows from \eqref{alow2}. The proof of
\eqref{acom2} is similar.

Due to \eqref{aep} the operators $z_{in}$ and $z_{ni}$ preserve 
the space 
$L(\lambda)^+$: for $i=1,\dots,n-1$
\begin{equation}\label{arlo}
z_{in}:L(\lambda)^+_{\mu}\to L(\lambda)^+_{\mu+\delta_i},\qquad
z_{ni}:L(\lambda)^+_{\mu}\to L(\lambda)^+_{\mu-\delta_i},
\end{equation}
where $\mu\pm\delta_i$ is obtained from $\mu$ by replacing $\mu_i$
with $\mu_i\pm 1$. 

The following is a key lemma in the derivation of 
the Gelfand--Tsetlin formulas.

\ble\label{lem:aximu}
Given $\mu$ satisfying
\eqref{amulam} the vector
\beq%\label{aximu}
\xi_{\mu}=z_{n1}^{\lambda^{}_1-\mu^{}_1}\cdots 
z_{n,n-1}^{\lambda^{}_{n-1}-\mu^{}_{n-1}} \ts \xi
\non
\end{equation}
spans the subspace $L(\lambda)^+_{\mu}$ and 
for each $i=1,\dots,n-1$ we have 
\beq\label{aact}
z_{in}\ts \xi_{\mu}=-(m_i-l_1)\cdots (m_i-l_n)\ts \xi_{\mu+\delta^{}_i},
\end{equation}
where
%\beq\label{aml}
$m_i=\mu_i-i+1$, $l_i=\lambda_i-i+1$.
%\end{equation}
It is supposed that $\xi_{\mu+\delta^{}_i}=0$ if $\lambda_i=\mu_i$.
\ele

\Proof By \eqref{arlo} the vector $\xi_{\mu}$ belongs to
the subspace $L(\lambda)^+_{\mu}$.
We need to show that it is nonzero. This will follow from
relations \eqref{aact}. We shall outline a proof 
of these relations which involves
the use of the Yangians; cf. \cite{z:gz}. 	
Consider the $n\times n$-matrix $E$ whose
$ij$-th entry is $E_{ij}$ and let $u$ be a formal variable.
Introduce the polynomial $T(u)$ with coefficients 
in the universal enveloping algebra $\U(\g_n)$:
\beq%\label{aqdet}
T(u)=\sum_{\sigma\in S_n}\sgn \sigma\ts (u+E)_{\sigma(1),1}\cdots
(u+E-n+1)_{\sigma(n),n}.
\non
\end{equation}
It is well known that all its coefficients belong
to the center of $\U(\g_n)$	(and generate the center); see 
e.g. \cite{hu:ci}.
This also easily follows from the properties the {\it quantum
determinant\/} of the Yangian for the Lie algebra $\gl(n)$;
see e.g. \cite{mno:yc}. 
Therefore, these coefficients act in $L(\lambda)$ as scalars
which can be easily found by applying $T(u)$ to the highest vector $\xi$:
\beq\label{atu}
T(u)|^{}_{L(\lambda)}=(u+l_1)\cdots (u+l_n).
\end{equation}
On the other hand, the center of $\U(\g_n)$ is
a subalgebra in the normalizer $\Norm \J$. 
We shall keep the same notation for the image of $T(u)$ in the
transvector algebra $\Z(\g_n,\g_{n-1})$.
To get explicit expressions of the coefficients of $T(u)$ in terms of
the raising and lowering operators we
consider $T(u)$ modulo the ideal $\J$ and apply the projection $p$.
The details can be found in 
\cite[Theorem~3.1]{m:yt}.
We have
\beq%\label{atzz}
T(u)=(u+E_{nn})\prod_{i=1}^{n-1}(u+h_i-1)-
\sum_{i=1}^{n-1}z_{in}z_{ni}\prod_{j=1,\ts j\ne i}^{n-1}
\frac{u+h_j-1}{h_i-h_j}.
\non
\end{equation}
In particular, $T(-h_i+1)=(-1)^{n-1} z_{in}z_{ni}$.

Now, \eqref{acom2} implies that $z_{in}\ts \xi_{\mu}=0$
unless $\lambda_i-\mu_i\geq 1$. In this case
using \eqref{acom}
we obtain
\beq%\label{acalc}
z_{in}\ts \xi_{\mu}=z_{in}z_{ni}\ts \xi_{\mu+\delta^{}_i}
=(-1)^{n-1}T(-h_i+1)\ts \xi_{\mu+\delta^{}_i}
= (-1)^{n-1} T(-m_i)\ts \xi_{\mu+\delta^{}_i},
\non
\end{equation}
where we have used $h_i\ts \xi_{\mu}=m_i\ts \xi_{\mu}$.
The relation \eqref{aact} now follows from \eqref{atu}.

Applying appropriate raising operators to the
vector $\xi_{\mu}$ we can obtain the highest vector $\xi$ of $L(\lambda)$
with a nonzero coefficient. This proves that $\xi_{\mu}\ne 0$.
\endproof

Given a Gelfand--Tsetlin pattern $\Lambda$ introduce the vector
$\xi^{}_{\Lambda}\in L(\lambda)$ by
\beq%\label{axiL}
\xi^{}_{\Lambda}=\prod_{k=2,\dots,n}^{\rightarrow}
\Bigl(z_{k1}^{\lambda^{}_{k1}-\lambda^{}_{k-1,1}}\cdots
z_{k,k-1}^{\lambda^{}_{k,k-1}-\lambda^{}_{k-1,k-1}}\Bigr)\ts\xi.
\non
\end{equation}

Lemma~\ref{lem:aximu} implies

\bco\label{cor:basis}
The vectors $\xi^{}_{\Lambda}$ parametrized by the patterns $\Lambda$
form a basis of the representation $L(\lambda)$.	\endproof
\eco

We now briefly outline
a derivation of formulas \eqref{aekk}--\eqref{aek+1k} 
which is standard;
see e.g. \cite{z:gz}. 
First, since $E_{nn}\ts z_{ni}=z_{ni}\ts (E_{nn}+1)$ for any $i$,
we have
\beq%\label{aenn}
E_{nn}\ts\xi_{\mu}=\Bigl(\sum_{i=1}^n\lambda_i-
\sum_{i=1}^{n-1}\mu_i\Bigr)\ts\xi_{\mu},
\non
\end{equation}
which implies \eqref{aekk}. To prove \eqref{aekk+1} is suffices
to find $E_{n-1,n}\ts \xi_{\mu\nu}$ where
\beq%\label{aen-1n}
\xi_{\mu\nu}=
z_{n-1,1}^{\mu^{}_1-\nu^{}_1}\cdots z_{n-1,n-2}^{\mu^{}_{n-2}-\nu^{}_{n-2}}
\ts\xi_{\mu},
\non
\end{equation}
and the $\nu_i$ satisfy
\beq%\label{anumu}
\mu_i-\nu_i\in\ZZ_+ 	\quad\text{and}\quad \nu_i-\mu_{i+1}\in
\ZZ_+ 	\qquad\text{for}\quad i=1,\dots,n-2.
\non
\end{equation}
Since $E_{n-1,n}$ commutes with the $z_{n-1,i}$, 
\beq%\label{aen-1n2}
E_{n-1,n}\ts \xi_{\mu\nu}=
z_{n-1,1}^{\mu^{}_1-\nu^{}_1}\cdots z_{n-1,n-2}^{\mu^{}_{n-2}-\nu^{}_{n-2}}\ts
E_{n-1,n}\ts\xi_{\mu}. 
\non
\end{equation}
Now use the following
identity in $\U'(\g_n)$ modulo the ideal $\J$ \cite{z:gz}
\beq\label{aen}
E_{n-1,n}=\sum_{i=1}^{n-1}z_{n-1,i}\ts z_{in}
\frac{1}{(h_i-h_1)\cdots\wedge\cdots(h_i-h_{n-1})},
\end{equation}
where $z_{n-1,n-1}:=1$. Applying \eqref{aact} we find that
\beq\label{aen-1n3}
E_{n-1,n}\ts \xi_{\mu\nu}=
-\sum_{i=1}^{n-1} \frac{(m_i-l_1)\cdots (m_i-l_n)}
{(m_i-m_1)\cdots\wedge\cdots (m_i-m_{n-1})}
\ts
\xi_{\mu+\delta_i,\nu}
\end{equation}
which proves \eqref{aekk+1}.
To prove \eqref{aek+1k} we use a contravariant 
bilinear form $\langle\ts,\rangle$
on $L(\lambda)$ uniquely determined by the
conditions:
\beq\label{aform}
\langle \xi,\xi\rangle=1,\qquad
\langle E_{ij}\ts\eta,\zeta\rangle=\langle \eta,E_{ji}\ts\zeta\rangle,
\quad \eta,\zeta\in L(\lambda).
\end{equation}
The basis $\{\xi_{\Lambda}\}$ is orthogonal with respect to this
form. This follows from \eqref{aact} and 
the fact that the operators $pE_{in}$ and $pE_{ni}$ are adjoint
to each other with respect to the restriction of the form $\langle\ts,\rangle$
to the space $L(\lambda)^+$. In particular, \eqref{aen-1n3} implies that
\beq%\label{aenn-1}
E_{n,n-1}\xi_{\mu\nu}=\sum_{i=1}^{n-1} c_i(\mu,\nu)\ts \xi_{\mu-\delta_i,\nu}
\non
\end{equation}
for some coefficients $c_i(\mu,\nu)$. Apply the operator $z_{j,n-1}$
to both sides of this relation. Since $z_{j,n-1}$ commutes with
$E_{n,n-1}$ we obtain from \eqref{aact} a recurrence relation
for the $c_i(\mu,\nu)$: if $\mu_j-\nu_j\geq 1$ then
\beq%\label{arec}
c_i(\mu,\nu+\delta_j)=c_i(\mu,\nu)\ts \frac{m_i-\gamma_j-1}{m_i-\gamma_j},
\non
\end{equation}
where $\gamma_j=\nu_j-j+1$.
This proves \eqref{aek+1k} by induction. \endproof

Note that
the original Gelfand--Tsetlin basis \cite{gt:fdu} is orthonormal.
The basis vectors in \cite{gt:fdu} 
coincide with the $\xi^{}_{\Lambda}$ up to
a norm factor which can be explicitly calculated; see e.g. \cite{z:gz}.

\section{A basis for odd orthogonal Lie algebras} \label{sec:boo}
\setcounter{equation}{0}

We shall enumerate the rows and columns of $(2n+1)\times (2n+1)$-matrices over
$\C$ by the indices $-n,\dots,-1,0,1,\dots,n$.

\subsection{Main theorem}\label{subsec:mt}

We keep the notation $E_{ij}$, $i,j=-n,\dots,n$ 
for the standard
basis of the Lie algebra $\gl(2n+1)$. 
Introduce the elements
\begin{equation}\label{Fij}
F_{ij}=E_{ij}-E_{-j,-i}. 
\end{equation}
We have $F_{-j,-i}=-F_{ij}$. In particular, $F_{-i,i}=0$ for all $i$.
The orthogonal Lie algebra $\g_n:=\oa(2n+1)$
can be identified
with the subalgebra in $\gl(2n+1)$
spanned by the elements $F_{ij}$, $i,j=-n,\dots,n$. 
The subalgebra $\g_{n-1}$ is spanned by the elements \eqref{Fij} with the
indices $i,j$ running over the set $\{-n+1,\dots,n-1\}$. 
Denote by $\h=\h_n$ the diagonal Cartan subalgebra in $\g_n$. 
The elements $F_{11},\dots,F_{nn}$ form a basis of $\h$. 

The finite-dimensional irreducible representations of $\g_n$
are in a one-to-one correspondence with $n$-tuples
$\lambda=(\lambda_1,\dots,\lambda_n)$ 
where all the entries $\lambda_i$ are simultaneously
integers or half-integers 
(elements of the set $\frac12+\ZZ$)
and the following
inequalities hold:
\begin{equation}%\label{inD}
0\geq\lambda_1\geq\lambda_2\geq\cdots\geq \lambda_n.
\non
\end{equation}
Such an $n$-tuple $\lambda$ is called the highest weight
of the corresponding representation which
we shall denote by $V(\lambda)$.
It contains a unique, up to a multiple, nonzero vector $\xi$
(the highest vector) such that
$
F_{ii}\ts\xi=\lambda_i\ts\xi
$
for $i=1,\dots,n$ and
$
F_{ij}\ts\xi=0
$
for $-n\leq i<j\leq n$.
Denote by $V(\lambda)^+$ the subspace of $\g_{n-1}$-highest vectors
in $V(\lambda)$:
\begin{equation}%\label{subhv}
V(\lambda)^+=\{\eta\in V(\lambda)\ |\ F_{ij}\ts \eta=0,
\qquad -n<i<j<n\}.
\non
\end{equation}
Given a $\g_{n-1}$-highest weight
$\mu=(\mu_1,\dots,\mu_{n-1})$ we denote by $V(\lambda)^+_{\mu}$
the corresponding weight subspace in $V(\lambda)^+$:
\begin{equation}%\label{hvwmu}
V(\lambda)^+_{\mu}=\{\eta\in V(\lambda)^+\ |\ F_{ii}\ts\eta=
\mu_i\ts\eta,\qquad i=1,\dots,n-1\}.
\non
\end{equation}
By the branching rule for the reduction $\g_n\downarrow\g_{n-1}$ \cite{z:cg}
we have
\beq\label{bra}
V(\lambda)|^{}_{\g_{n-1}}\simeq\bigoplus c(\mu)\ts V'(\mu),
\end{equation}
where $V'(\mu)$ is the irreducible finite-dimensional representation of
$\g_{n-1}$ with the highest weight $\mu$, and $c(\mu)$ equals the number
of $n$-tuples $\rho=(\rho_1,\dots,\rho_n)$ satisfying the inequalities
\beq\label{ineq}
\begin{split}
&-\lambda_1\geq\rho_1\geq\lambda_1\geq\rho_2\geq\lambda_2\geq \cdots\geq
\rho_{n-1}\geq\lambda_{n-1}\geq\rho_n\geq\lambda_n,\\
&-\mu_1\geq\rho_1\geq\mu_1\geq\rho_2\geq\mu_2\geq \cdots\geq
\rho_{n-1}\geq\mu_{n-1}\geq\rho_n
\end{split}
\end{equation}
with all the $\rho_i$ and $\mu_i$ being simultaneously
integers or half-integers together with the $\lambda_i$.
Any nonzero vector in $V(\lambda)^+_{\mu}$ generates a
$\g_{n-1}$-submodule in $V(\lambda)$ isomorphic to $V'(\mu)$.
We obviously have $\dim V(\lambda)^+_{\mu}=c(\mu)$.
Basis vectors in $V(\lambda)^+_{\mu}$ can be
parametrized by the $n$-tuples $\rho$. 
We shall be using an equivalent parameterization
by $(n+1)$-tuples $\nu=(\sigma,\nu_1,\dots,\nu_n)$,
where $\nu_i=\rho_i$ for $i\geq 2$, and 
\beq
(\sigma,\nu_1)=\begin{cases}
(0,\rho_1)\quad&\text{if $\rho_1\leq 0$},\\
(1,-\rho_1)\quad&\text{if $\rho_1> 0$}.
\end{cases}
\non
\end{equation}
A parameterization of basis vectors in $V(\lambda)$ is obtained by
using its subsequent restrictions to the subalgebras of the chain
$
\g_1\subset\g_2\subset\cdots\subset\g_{n-1}\subset\g_n.
$
Define a {\it %Gelfand--Tsetlin 
pattern\/} $\Lambda$ associated with
$\lambda$ as an array of the form
\begin{align}
\sigma^{}_{n}\quad&\qquad\lambda^{}_{n1}\qquad\lambda^{}_{n2}
\qquad\qquad\cdots\qquad\qquad\lambda^{}_{nn}\non\\
&\lambda'_{n1}\qquad \lambda'_{n2}
\qquad\qquad\cdots\qquad\qquad\lambda'_{nn}\non\\
\sigma^{}_{n-1}\ &\qquad\lambda^{}_{n-1,1}\qquad\cdots
\qquad\lambda^{}_{n-1,n-1}\non\\
&\lambda'_{n-1,1}
\qquad\cdots\qquad\lambda'_{n-1,n-1}\non\\
\qquad\cdots&\qquad\cdots\qquad\cdots\non\\
\sigma^{}_{1}\quad&\qquad\lambda^{}_{11}\non\\
&\lambda'_{11}\non
\end{align}
such that $\lambda=(\lambda^{}_{n1},\dots, \lambda^{}_{nn})$,
each $\sigma^{}_{k}$ is $0$ or $1$, the remaining
entries are all non-positive
integers or non-positive half-integers together with the
$\lambda_i$, and the following inequalities hold
$$
\lambda'_{k1}\geq\lambda^{}_{k1}\geq\lambda'_{k2}\geq
\lambda^{}_{k2}\geq \cdots\geq
\lambda'_{k,k-1}\geq\lambda^{}_{k,k-1}\geq\lambda'_{kk}\geq\lambda^{}_{kk}
$$
for $k=1,\dots,n$; and
$$
\lambda'_{k1}\geq\lambda^{}_{k-1,1}\geq\lambda'_{k2}\geq
\lambda^{}_{k-1,2}\geq \cdots\geq
\lambda'_{k,k-1}\geq\lambda^{}_{k-1,k-1}\geq\lambda'_{kk}
$$
for $k=2,\dots,n$; in addition, in the case of integer $\lambda_i$
the condition
$$
\lambda'_{k1}\leq -1\qquad\text{if}\quad \sigma^{}_{k}=1
$$
should hold for all $k=1,\dots,n$.
Let us set $l^{}_{k0}=-1/2$ for all $k$, and
$$
l^{}_{ki}=\lambda^{}_{ki}-i+1/2,\qquad l'_{ki}=\lambda'_{ki}-i+1/2,\qquad
1\leq i\leq k\leq n.
$$
Given a pattern $\Lambda$ set for $i=0,1,\dots,k-1$
\beq
A_{ki}=\prod_{a=1, \ts a\ne i}^{k-1}
\frac{1}{l^{}_{k-1,i}-l^{}_{k-1,a}}\cdot
\prod_{a=1}^{k-1}
\frac{1}{l^{}_{k-1,i}+l^{}_{k-1,a}}.
\non
\end{equation}
Furthermore,
introduce polynomials $B_{ki}(x)$ by
\beq
B_{ki}(x)=\prod_{a=1, \ts a\ne i}^{k}
\frac{(x+l'_{ka}+1)(x-l'_{ka})}{l'_{ka}-l'_{ki}},
\non
\end{equation}
and define the numbers $C_{ki}$ by
\beq
C_{ki}=
l'_{ki}\ts(1-2\ts\sigma^{}_k-2\ts l'_{ki})\ts
\prod_{a=1}^k (l^{}_{ka}-l'_{ki})
\prod_{a=1}^{k-1} (l^{}_{k-1,a}-l'_{ki})
\prod_{a=1,\ts a\ne i}^{k} \frac{1}{l'_{ka}-l'_{ki}}.\non
\end{equation}
We denote by
$\Lambda\pm\delta^{}_{ki}$ and $\Lambda+\delta'_{ki}$ the arrays
obtained from $\Lambda$ by replacing $\lambda^{}_{ki}$ and $\lambda'_{ki}$
by $\lambda^{}_{ki}\pm1$ and $\lambda'_{ki}+1$ respectively.

The following is our main theorem which will be proved
in Sections~\ref{subsec:cb} and
\ref{subsec:me}.

\bth\label{thm:main} There exists a basis
$\{\zeta^{}_{\Lambda}\}$ of $V(\lambda)$
parametrized by the patterns $\Lambda$
such that the action of the generators
of $\g_n$
is given by the formulas
\begin{align}
F^{}_{kk}\ts \zeta^{}_{\Lambda}&=\left(\sigma^{}_{k}+
2\sum_{i=1}^{k}\lambda'_{ki}-
\sum_{i=1}^k\lambda^{}_{ki}-\sum_{i=1}^{k-1}\lambda^{}_{k-1,i}\right)
\zeta^{}_{\Lambda},
\non\\
F^{}_{k-1,-k}\ts \zeta^{}_{\Lambda}&=A_{k0}\ts\zeta^{}_{\Lambda}(k,0)
+\sum_{i=1}^{k-1}
A_{ki}\left(\frac{1}{l^{}_{k-1,i}+1/2}\ts\zeta^+_{\Lambda}(k,i)-
\frac{1}{l^{}_{k-1,i}-1/2}\ts\zeta^-_{\Lambda}(k,i)\right).
\non
\end{align}
Here the following notation has been used
\begin{align}
\zeta^-_{\Lambda}(k,i)&=
\zeta^{}_{\Lambda-\delta^{}_{k-1,i}},\non\\
\zeta^+_{\Lambda}(k,i)&=\sum_{j=1}^{k}\sum_{m=1}^{k-1}
B_{kj}(l_{k-1,i})B_{k-1,m}(l_{k-1,i})\ts
\zeta^{}_{\Lambda+\delta'_{kj}+\delta^{}_{k-1,i}+\delta'_{k-1,m}},
\non
\end{align}
and $\zeta^{}_{\Lambda}:=0$ if $\Lambda$ is not a pattern.
Furthermore,
\begin{alignat}{2}
\zeta^{}_{\Lambda}(k,0)&=(-1)^k\ts \zeta^{}_{\bar{\Lambda}}
\quad&&\text{if}\quad
\sigma^{}_{k}=\sigma^{}_{k-1}=0,\non\\
{}&=\sum_{j=1}^{k}B_{kj}(l_{k-1,0})\ts
\zeta^{}_{\bar{\Lambda}+\delta'_{kj}}
\quad&&\text{if}\quad
\sigma^{}_{k}=1,\ \sigma^{}_{k-1}=0,\non\\
{}&=-\sum_{m=1}^{k-1}
B_{k-1,m}(l_{k-1,0})\ts
\zeta^{}_{\bar{\Lambda}+\delta'_{k-1,m}}
\quad&&\text{if}\quad
\sigma^{}_{k}=0,\ \sigma^{}_{k-1}=1,\non\\
{}&=(-1)^{k-1}\sum_{j=1}^{k}\sum_{m=1}^{k-1}
B_{kj}(l_{k-1,0})B_{k-1,m}(l_{k-1,0})\ts&&
\zeta^{}_{\bar{\Lambda}+\delta'_{kj}+\delta'_{k-1,m}}\non\\
&\quad&&\text{if}\quad
\sigma^{}_{k}=\sigma^{}_{k-1}=1,
\non
\end{alignat}
where
$\bar{\Lambda}$ is obtained from $\Lambda$ by replacing
$\sigma^{}_{k}$ and $\sigma^{}_{k-1}$ respectively with
$\sigma^{}_{k}+1$ and $\sigma^{}_{k-1}+1$ {\rm(}modulo $2${\rm)}.
The action of $F^{}_{k-1,k}$ is found from the relation
\begin{equation}
F^{}_{k-1,k}=\Bigl[\Phi_{k-1,-k}(u+2)\ts\Phi_{-k,k}-
\Phi_{-k,k}\Phi_{k-1,-k}(u)\Bigr]^{}_{u=0},\non
\end{equation}
where the operator $\Phi_{-k,k}$ acts on the basis elements by the rule
\begin{equation}
\Phi_{-k,k}\ts\zeta^{}_{\Lambda}=\sum_{i=1}^{k} C_{ki}\ts
(F^{}_{kk}-l'_{ki}+1)\ts\zeta^{}_{\Lambda-\delta'_{ki}}\non
\end{equation}
while the action of $\Phi_{k-1,-k}(u)$ is given by
\begin{multline}
\Phi_{k-1,-k}(u)\ts \zeta^{}_{\Lambda}=
\frac{A_{k0}}{u+F^{}_{kk}-3/2}\ts \zeta^{}_{\Lambda}(k,0)\\
+\sum_{i=1}^{k-1}
A_{ki}\Bigg(\frac{1}{(l^{}_{k-1,i}+1/2)(u+l^{}_{k-1,i}+F^{}_{kk}-1)}
\ts\zeta^+_{\Lambda}(k,i)\\
-\frac{1}{(l^{}_{k-1,i}-1/2)(u-l^{}_{k-1,i}+F^{}_{kk}-1)}
\ts\zeta^-_{\Lambda}(k,i)\Bigg).
\non
\end{multline}
\eth

\medskip
\noindent
{\it Remark\/}.  The image of $\zeta^{}_{\Lambda}$ under
the operator $\Phi_{k-1,-k}(u+2)\ts\Phi_{-k,k}-
\Phi_{-k,k}\Phi_{k-1,-k}(u)$ at $u=0$ may be undefined
for some patterns $\Lambda$.
To get the action of $F^{}_{k-1,k}$,
one should first calculate its matrix elements in a ``generic"
representation $V(\lambda)$ and then specialize the parameters;
see Section~\ref{subsec:me}. For example,
consider the case $n=1$. It will be shown
in Section~\ref{subsec:cb} that the basis vectors 
in $V(\lambda)$ are given by
\beq
\zeta^{}_{\Lambda}=F^{\sigma_1}_{10}\ts 
(F^{}_{10}F^{}_{0,-1})^{\lambda'_{11}-\lambda^{}_{11}}_{}\ts\xi.
\non
\end{equation}
Furthermore, the operators $\Phi_{-1,1}$ and $\Phi_{0,-1}(u)$
are defined by
\begin{align}
\Phi_{-1,1}&=-\frac12\ts F_{01}^2,
\non\\
\Phi_{0,-1}(u)&=F^{}_{0,-1}\ts\frac{1}{u+F_{11}-1/2};
\non
\end{align}
see Section~\ref{subsec:me}. In the case $\lambda=(-1/2)$
the basis of  $V(\lambda)$ consists of two vectors
$\xi$ and $\xi'=F_{10}\ts\xi$. Therefore, $\Phi_{-1,1}$
is the zero operator in $V(\lambda)$ while the image
$\Phi_{0,-1}(0)\xi'$ is not defined. On the other hand,
we find directly that $F^{}_{01}\ts\xi'=1/2\ts \xi$.
\endproof

\subsection{Transvector algebra $\Z(\g_n,\g_{n-1})$}

Consider the extension 
of the universal enveloping algebra $\U(\g_n)$
\begin{equation}%\label{ext}
\U'(\g_n)=\U(\g_n)\ot_{\U(\h)} \R(\h),
\non
\end{equation}
where $\R(\h)$ is the field of fractions of the commutative algebra
$\U(\h)$. Let $\J$
denote the left ideal in $\U'(\g_n)$ generated by 
the elements $F_{ij}$ with $-n<i<j<n$.
The {\it transvector algebra\/} $\Z(\g_n,\g_{n-1})$ is the quotient algebra
of the normalizer
$$
\Norm \J=\{x\in\U'(\g_n)\ |\ \J x\subseteq \J\}
$$
modulo the two-sided ideal $\J$ \cite{z:it}.
It is an algebra over $\C$ and an $\R(\h)$-bimodule.
Let $p=p_{n-1}$ denote the {\it extremal projection\/}
for the Lie algebra $\g_{n-1}$ \cite{ast:po,z:it}; cf. Section~\ref{sec:gtb}.
It satisfies the following (characteristic) relations
\beq%\label{pchar}
F_{ij}\ts p=p\ts F_{ji}=0\qquad\text{for}\quad -n<i<j<n. 
\non
\end{equation}
The projection $p$ naturally acts in the space
$\U'(\g_n)/\J$ and its image coincides with $\Z(\g_n,\g_{n-1})$.
The elements
\begin{equation}\label{gener}
pF_{ia}=-pF_{-a,-i},\qquad a=-n,n,\quad i=-n+1,\dots,n-1
\end{equation}
are generators of $\Z(\g_n,\g_{n-1})$ \cite{z:it}. 
Set
\begin{equation}%\label{fiD}
f_0=-1/2,\qquad f_i=F_{ii}-i+1/2,\qquad f_{-i}=-f_i
\non
\end{equation}
for $i=1,\dots,n$; and set $f'_i=-f_{-i}$ for all $i$.
The elements \eqref{gener} can be given by
the following explicit formulas (modulo $\J$):
\begin{equation}%\label{pFia}
pF_{ia}=\sum_{i>i_1>\cdots>i_s>-n}
F_{ii_1}F_{i_1i_2}\cdots F_{i_{s-1}i_s}F_{i_sa}
\frac{1}{(f_i-f_{i_1})\cdots (f_i-f_{i_s})},
\non
%\\ \label{pFai} pF_{ai}&=\sum_{i<i_1<\cdots<i_s<n}
%F_{i_1i}F_{i_2i_1}\cdots F_{i_si_{s-1}}F_{ai_s}
%\frac{1}{(f_i-f_{i_1})\cdots (f_i-f_{i_s})},\non
\end{equation}
where $s=0,1,\dots$. We shall use normalized
generators of $\Z(\g_n,\g_{n-1})$ defined by
\begin{align}\label{zia}
z_{ia}&=pF_{ia}\ts (f_i-f_{i-1})\cdots(f_i-f_{-n+1}),\\
\label{zai}
z_{ai}&=pF_{ai}\ts (f'_i-f'_{i+1})\cdots(f'_i-f'_{n-1}).
\end{align}
We obviously have $z_{ai}=(-1)^{n-i}\ts z_{-i,-a}$.
The elements $z_{ia}$ satisfy 
certain quadratic relations \cite{z:it}. We shall use
the following ones below (cf. \eqref{acom} and \eqref{acom2}):
for $a,b\in\{-n,n\}$ and $i+j\ne 0$ one has
\begin{equation}\label{relcom}
z_{aj}z_{bi}(f'_i-f'_j+1)=z_{bi}z_{aj}(f'_i-f'_j)
+z_{ai}z_{bj}. 
%\non
\end{equation}
In particular, $z_{ai}$ and $z_{aj}$ commute for $i+j\ne 0$.
One easily verifies that $z_{ai}$ and $z_{bi}$ also commute for
$i\ne 0$ and all $a,b$.

%The following equivalent formula holds for the $z_{ai}$:
%\beq\label{zaieq}
%z_{ai}=(f_i-f_{j_1})\cdots(f_i-f_{j_k})
%\sum_{n>i_1>\cdots>i_s>i}
%F_{ai_1}F_{i_1i_2}\cdots F_{i_{s-1}i_s}F_{i_si},
%\end{equation}
%where $\{j_1,\dots,j_k\}$ is the complement to
%$\{i_1,\dots,i_s\}$ in the set $\{i+1,\dots,n-1\}$.
The elements $z_{ia}$ and $z_{ai}$ 
naturally act in the space 
$V(\lambda)^+$ and are called
the {\it raising\/} and {\it lowering operators\/}.
One has for $i=1,\dots,n-1$:
\begin{equation}%\label{rlo}
z_{ia}:V(\lambda)^+_{\mu}\to V(\lambda)^+_{\mu+\delta_i},\qquad
z_{ai}:V(\lambda)^+_{\mu}\to V(\lambda)^+_{\mu-\delta_i},\non
\end{equation}
where $\mu\pm\delta_i$ is obtained from $\mu$ by replacing $\mu_i$
with $\mu_i\pm 1$. 
The operators $z^{}_{0a}$ preserve each
subspace $V(\lambda)^+_{\mu}$. 

We shall need the following element which can be checked to belong
to the normalizer $\Norm \J$, and so it can be regarded as
an element of the algebra $\Z(\g_n,\g_{n-1})$:
\begin{equation}\label{zn-n}
z_{n,-n}=
\sum_{n>i_1>\cdots>i_s>-n}
F_{ni_1}F_{i_1i_2}\cdots F_{i_s,-n}\ts
(f_{n}-f_{j_1})\cdots (f_{n}-f_{j_k}),
\end{equation}
where $s=1,2,\dots$ and $\{j_1,\dots,j_k\}$ is the complement
to the subset $\{i_1,\dots,i_s\}$ in
$\{-n+1,\dots,n-1\}$. 
The following relation is proved exactly as its
$C$ and $D$ series counterparts \cite{m:br, m:wb}
(cf. \eqref{aen}): for $a=-n,n$
\begin{equation}\label{Fn-1a}
F_{n-1,a}=\sum_{i=-n+1}^{n-1}z_{n-1,i}\ts z_{ia}
\frac{1}{(f_i-f_{-n+1})\cdots\wedge\cdots(f_i-f_{n-1})},
\end{equation}
where $z_{n-1,n-1}:=1$ and the equalities are considered 
in $\U'(\g_n)$ modulo
the ideal $\J$.

\subsection{Yangians and twisted Yangians}

Let us introduce the $\gl(2)$-{\it Yangian\/} $\Y(2)$ and the 
({\it orthogonal\/})
{\it twisted Yangian\/} $\Y^+(2)$; see \cite{mno:yc} for more details.
The Yangian $\Y(2)$ is the
complex associative algebra with the
generators $t_{ab}^{(1)},t_{ab}^{(2)},\dots$ where 
$a,b\in\{-n,n\}$,
and the defining relations
\begin{equation}\label{rel}
[t_{ab}(u),t_{cd}(v)]=\frac{1}{u-v}
\Big(t_{cb}(u)t_{ad}(v)-t_{cb}(v)t_{ad}(u)\Big),
\end{equation}
where
\begin{equation}\label{ser}
t_{ab}(u): = \delta_{ab} + t^{(1)}_{ab} u^{-1} + t^{(2)}_{ab}u^{-2} +
\cdots \in \Y(2)[[u^{-1}]].
%\non
\end{equation}
Introduce the series $s_{ab}(u)$, $a,b\in\{-n,n\}$ by
\begin{equation}%\label{sab}
s_{ab}(u)=t_{an}(u)t_{-b,-n}(-u)+t_{a,-n}(u)t_{-b,n}(-u).
\non
\end{equation}
Write
$
s_{ab}(u)=\delta_{ab}+s_{ab}^{(1)}u^{-1}+s_{ab}^{(2)}u^{-2}+\cdots.
$
The twisted Yangian $\Y^+(2)$ is defined as the subalgebra of $\Y(2)$
generated by the elements $s_{ab}^{(1)},s_{ab}^{(2)},\dots$ where 
$a,b\in\{-n,n\}$. Also, $\Y^+(2)$ can be viewed as
an abstract algebra with generators $s_{ab}^{(r)}$ and the following
defining relations (see \cite[Section~3]{mno:yc}):
\begin{align}%\label{deftw}
[s^{}_{ab}(u),s^{}_{cd}(v)]&={1\over u-v}
\Big(s^{}_{cb}(u)s^{}_{ad}(v)-s^{}_{cb}(v)s^{}_{ad}(u)\Big)\non\\
\label{deftw}
{}&-{1\over u+v}\Big(s^{}_{a,-c}(u)s^{}_{-b,d}(v)-
s^{}_{c,-a}(v)s^{}_{-d,b}(u)\Big)\\
{}&+{1\over u^2-v^2}\Big(s^{}_{c,-a}(u)s^{}_{-b,d}(v)-
s^{}_{c,-a}(v)s^{}_{-b,d}(u)\Big) 
\non
\end{align}
and
\beq\label{sym}
s^{}_{-b,-a}(-u)=\frac{2u+1}{2u}\ts
s^{}_{ab}(u)-\frac{1}{2u}\ts s^{}_{ab}(-u). 
\end{equation}
The Yangian $\Y(2)$ is a Hopf algebra with the 
coproduct 
\begin{equation}\label{cop}
\Delta (t_{ab}(u))=t_{an}(u)\ot
t_{nb}(u)+t_{a,-n}(u)\ot t_{-n,b}(u).
\end{equation}
The twisted Yangian $\Y^+(2)$ is a left coideal in $\Y(2)$
with
\begin{equation}\label{cops}
\Delta (s_{ab}(u))=\sum_{c,d\in\{-n,n\}}t_{ac}(u)t_{-b,-d}(-u)\ot
s_{cd}(u).
\end{equation}

Given a pair of complex numbers $(\alpha,\beta)$ 
such that $\alpha-\beta\in\ZZ_+$
we denote by 
$L(\alpha,\beta)$ the irreducible representation of the Lie algebra
$\gl(2)$ with the highest weight $(\alpha,\beta)$ with respect to the
upper triangular Borel subalgebra; see Section~\ref{sec:gtb}. 
We have 
$\dim L(\alpha,\beta)=\alpha-\beta+1$. We may regard $L(\alpha,\beta)$ as
a $\Y(2)$-module by using the algebra homomorphism $\Y(2)\to\U(\gl(2))$
given by
\begin{equation}%\label{hom}
t_{ab}(u)\mapsto \delta_{ab}+E_{ab}u^{-1},\qquad a,b\in\{-n,n\}.
\non
\end{equation}
The coproduct \eqref{cop} allows one to construct representations of
$\Y(2)$ of the form
\begin{equation}%\label{tenpr}
L=L(\alpha_1,\beta_1)\ot\cdots\ot L(\alpha_n,\beta_n).
\non
\end{equation}
Note that the generators $t_{ab}^{(r)}$ with $r>n$ act as zero operators
in $L$. Therefore, the operators $T_{ab}(u)=u^n\ts t_{ab}(u)$
are polynomials in $u$:
\begin{equation}%\label{Tab}
T_{ab}(u)=\delta_{ab}u^n+t_{ab}^{(1)}u^{n-1}+\cdots+t_{ab}^{(n)}.
\non
\end{equation}

For any $\gamma\in\C$ denote by $W(\gamma)$ the one-dimensional
representation of $\Y^+(2)$ spanned by a vector $w$ such that
\begin{equation}%%\label{onedim}
s_{nn}(u)\ts w=\frac{u+\gamma}{u+1/2}\ts w,\qquad 
s_{-n,-n}(u)\ts w=\frac{u-\gamma+1}{u+1/2}\ts w,
\non
\end{equation}
and $s_{a,-a}(u)\ts w=0$ for $a=-n,n$. 
By \eqref{cops} we can regard the tensor
product $L\ot W(\gamma)$ as a representation of $\Y^+(2)$.
Representations of this type
essentially exhaust all finite-dimensional irreducible representations of
$\Y^+(2)$ \cite{m:fd}. The vector space isomorphism
\begin{equation}\label{isomL}
L\ot W(\gamma)\to L,\qquad v\ot w\mapsto v,\qquad v\in L
\end{equation}
provides $L\ot W(\gamma)$ with an action of $\Y(2)$.

\subsection{Construction of the basis}\label{subsec:cb}

Introduce the following polynomials in $u$ with coefficients in the
transvector algebra $\Z(\g_n,\g_{n-1})$: for $a,b\in\{-n,n\}$
%%% 
%%%  The elements $Z_{ab}(u)$ are $\tilde Z_{ab}(u)$
%%%  from the notebook.
%%%
\beq\label{Zab}
Z_{ab}(u)=-\Big(\delta_{ab}(u-n+1)+F_{ab}\Big)
\prod_{i=-n+1}^{n-1}(u+g_i)
+\sum_{i=-n+1}^{n-1}z_{ai}z_{ib}
\prod_{j=-n+1,\ts j\ne  i}^{n-1}\frac{u+g_j}{g_i-g_j},
\end{equation}
where $g_i:=f_i+1/2$ for all $i$.

\bpr\label{prop:map}
The mapping
\begin{equation}\label{map}
s_{ab}(u)\mapsto -u^{-2n}\ts Z_{ab}(u),\qquad a,b\in\{-n,n\}
\end{equation}
defines an algebra homomorphism $\Y^+(2)\to \Z(\g_n,\g_{n-1})$.
\epr

\Proof One of the possible ways to prove the claim is to
check directly that the relations \eqref{deftw} and \eqref{sym}
are satisfied with the $s_{ab}(u)$ respectively replaced with
$Z_{ab}(u)$. Here one needs to use the quadratic relations
in the transvector algebra $\Z(\g_n,\g_{n-1})$.
In addition to \eqref{relcom} the relations which express
$z_{ia}z_{bi}$ in terms of the $z_{bj}z_{ja}$ are needed;
see \cite{z:it}.

Alternatively, we can follow the approach of \cite[Section~5]{m:br}
to construct first a homomorphism from $\Y^+(2)$ to the centralizer
$\CC_n$ of $\g_{n-1}$ in $\U(\g_n)$ and then calculate
the images of the centralizer elements in the
algebra  $\Z(\g_n,\g_{n-1})$. 
The calculation is similar to that
in the symplectic case \cite{m:br}; see also \cite{m:wb}.
We shall only give a few key formulas here.
Introduce the
$(2n+1)\times (2n+1)$-matrix $F=(F_{ij})$ whose $ij$th entry is
the element $F_{ij}\in\g_n$ and set
\begin{equation}%\label{Fu}
 F(u)=1+\frac{F}{u+1/2}.
\non
\end{equation}
Denote by $\wh { F}(u)$ the corresponding {\it Sklyanin comatrix\/};
see \cite[Section~2]{m:fd}. The mapping
\begin{equation}\label{homom}
s_{ab}(u)\mapsto  c(u)\ts  \wh { F}(-u+n-1/2)_{ab},\qquad a,b\in\{-n,n\},
\end{equation}
where $c(u)=(1-u^{-2})(1-4\ts u^{-2})\cdots(1-(n-1)^2\ts u^{-2})$,
defines an algebra homomorphism $\Y^+(2)\to\CC_n$ 
\cite[Proposition~2.1]{m:fd};
cf. \cite{o:ty}. Its composition with the natural
homomorphism $\CC_n\to \Z(\g_n,\g_{n-1})$ gives \eqref{map}.
\endproof

As it follows from the branching rule \eqref{bra},
%(Corollary~\ref{cor:bra})
the space $V(\lambda)^+_{\mu}$ is nonzero
only if
there exists $\nu$ such that
the inequalities 
\eqref{ineq} hold.
We shall be assuming that this condition is satisfied.
Proposition~\ref{prop:map} allows one to equip
$V(\lambda)^+_{\mu}$ with a structure of a $\Y^+(2)$-module
defined via the homomorphism \eqref{map}. 
The next theorem provides an identification of this 
module.

\bth\label{thm:isom}
The
$\ts\Y^+(2)$-module $V(\lambda)^+_{\mu}$ is isomorphic
to the direct sum of two irreducible submodules,
%\begin{equation}\label{isom}
$V(\lambda)^+_{\mu}\simeq U\oplus U'$,
%\end{equation}
where
\begin{align}\label{uu}
U&=L(0,\beta_1)\ot L(\alpha_2,\beta_2)\ot\cdots\ot 
L(\alpha_{n},\beta_{n})\ot W(1/2),\\
%\label{uus}
U'&=L(-1,\beta_1)\ot L(\alpha_2,\beta_2)\ot\cdots\ot 
L(\alpha_{n},\beta_{n})\ot W(1/2),
\end{align}
if the $\lambda_i$ are integers 
{\rm(}it is supposed that $U'=\{0\}$ if $\beta_1=0${\rm)}; or
\begin{align}%\label{uuws2}
U&=L(-1/2,\beta_1)\ot L(\alpha_2,\beta_2)\ot\cdots\ot 
L(\alpha_{n},\beta_{n})\ot W(0),
\\
\label{uu2}
U'&=L(-1/2,\beta_1)\ot L(\alpha_2,\beta_2)\ot\cdots\ot 
L(\alpha_{n},\beta_{n})\ot W(1),
\end{align}
if the $\lambda_i$ are half-integers, and the following notation
is used
\begin{alignat}{2}%\label{alphai}
\alpha_i&=\min\{\lambda_{i-1},\mu_{i-1}\}-i+1, \qquad&&i=2,\dots,n,
\non
\\
\beta_i&=\max\{\lambda_{i},\mu_{i}\}-i+1, \qquad&&i=1,\dots,n,
\non
\end{alignat}
with $\mu_{n}:=-\infty$.
In particular, each of $U$ and $U'$ {\rm(}and hence 
$V(\lambda)^+_{\mu}${\rm)} is equipped with
an action of $\Y(2)$ defined by \eqref{isomL}.
\eth

\Proof Consider the following two vectors in $V(\lambda)^+_{\mu}$
\begin{equation}\label{hv}
\xi_{\mu}=\prod_{i=1}^{n-1}\Big(z_{ni}^{\max\{\lambda_i,\mu_i\}-\mu_i}
z_{i,-n}^{\max\{\lambda_i,\mu_i\}-\lambda_i}\Big)\ts \xi,\qquad
\xi'_{\mu}=z^{}_{n0}\ts \xi_{\mu}.
\end{equation}
Repeating the arguments of the proof of Theorem~5.2 in \cite{m:br}
we can show that both $\xi_{\mu}$ and $\xi'_{\mu}$
are eigenvectors for $s_{nn}(u)$ and are annihilated by
$s_{-n,n}(u)$.
Namely, 
\beq\label{eig}
s_{nn}(u)\xi_{\mu}=\mu(u)\xi_{\mu},\qquad 
s_{nn}(u)\xi'_{\mu}=(1+u^{-1})\mu(u)\xi'_{\mu},
\end{equation}
where
\begin{equation}%\label{muu}
\mu(u)=(1-\alpha_2u^{-1})\cdots(1-\alpha_{n}u^{-1})
(1+\beta_1u^{-1})\cdots
(1+\beta_{n}u^{-1}).
\non
\end{equation}
This is proved simultaneously with the following relations
by induction on the degree of the monomial in \eqref{hv}:
for $i=1,\dots,n-1$
\begin{equation}\label{zinact}
z_{in}\ts\xi_{\mu}=-(m_i+\alpha_1)\cdots (m_i+\alpha_n)
(m_i-\beta_1)\cdots (m_i-\beta_{n})\ts \xi_{\mu+\delta_i},
%\non
\end{equation}
and
\begin{equation}\label{z-niact}
z_{-ni}\ts\xi_{\mu}=-(m_i-\alpha_1-1)\cdots(m_i-\alpha_{n}-1)
(m_i+\beta_1-1)\cdots 
(m_i+\beta_{n}-1)\ts \xi_{\mu-\delta_i},
%\non
\end{equation}
where $\alpha_1=0$ and
$m_i=\mu_i-i+1$ for $i=1,\dots,n-1$.
Indeed, we note that if $\mu_i\geq \lambda_i$ then $z_{in}\ts\xi_{\mu}=0$
which is implied by \eqref{relcom}. This agrees with
\eqref{zinact} because in this case $\beta_i=m_i$. Now assume that
$\mu_i<\lambda_i$. We have
%\beq%\label{zZ}
$z_{in}\ts\xi_{\mu}=z_{in}z_{ni}\ts\xi_{\mu+\delta_i}$
%\non
%\end{equation}
by \eqref{relcom}. Formula \eqref{Zab} gives 
%\beq%\label{zzZ}
$z_{in}z_{ni}=z_{-n,-i}z_{-i,-n}=Z_{-n,-n}(-g_{-i})$.
%\non
%\end{equation}
Further,
\beq%\label{Zxi}
Z_{-n,-n}(-g_{-i})\ts\xi_{\mu+\delta_i}=
Z_{-n,-n}(m_i)\ts\xi_{\mu+\delta_i}.
\non
\end{equation}
By Proposition~\ref{prop:map} and the symmetry relation \eqref{sym}
we can write
\beq%\label{symZ}
Z_{-n,-n}(m_i)=\frac{2m_i-1}{2m_i}Z_{nn}(-m_i)+
\frac{1}{2m_i}Z_{nn}(m_i).
\non
\end{equation}
By induction, $Z_{nn}(u)\ts \xi_{\mu+\delta_i}$ can be found
from \eqref{map} and \eqref{eig} which gives \eqref{zinact}. The proof of
\eqref{z-niact} is very similar. To prove \eqref{eig} we apply
the induction hypotheses to \eqref{zinact} and \eqref{z-niact}
and also use the relation
\beq\label{zon} 
z^{}_{0n}\ts\xi_{\mu}=0
\end{equation} 
which is
a consequence of \eqref{relcom}. The relations
\beq\label{ZZxi}
Z_{-n,n}(u)\ts\xi_{\mu}=0,\qquad Z_{-n,n}(u)\ts\xi'_{\mu}=0
\end{equation} 
follow from \eqref{Zab},  \eqref{zinact} and \eqref{z-niact}.

Both vectors $\xi_{\mu}$ and $\xi'_{\mu}$ are nonzero, except for the case
$\beta_1=0$ where $\xi'_{\mu}=0$. Indeed, applying
appropriate operators $z_{in}$ to $\xi_{\mu}$ or $\xi'_{\mu}$ repeatedly,
we can obtain the highest vector $\xi$ of $V(\lambda)$ 
with a nonzero coefficient.
It follows from \cite[Corollary~6.6]{m:fd} that
the tensor products \eqref{uu}--\eqref{uu2}
are irreducible representations
of $\Y^+(2)$. An easy calculation shows that
the highest weights of the $\Y^+(2)$-modules $U$ and $U'$ 
respectively coincide
with the $\Y^+(2)$-weights of the vectors
$\xi_{\mu}$ and $\xi'_{\mu}$. So, $U$ and $U'$ are 
respectively isomorphic
to quotients of the $\Y^+(2)$-submodules in $V(\lambda)^+_{\mu}$
generated by $\xi_{\mu}$ and $\xi'_{\mu}$.
On the other hand, the branching rule \eqref{bra} implies that
\beq%\label{dim}
\dim V(\lambda)^+_{\mu}=\dim U+\dim U'.
\non
\end{equation}
Therefore, to complete the proof of the theorem we need to show that
the $\Y^+(2)$-submodules generated by  $\xi_{\mu}$ and $\xi'_{\mu}$
are disjoint. For this we employ 
a contravariant 
bilinear form $\langle\ts,\rangle$
on $V(\lambda)$ uniquely determined by the
conditions:
\beq%\label{form}
\langle \xi,\xi\rangle=1,\qquad
\langle F_{ij}\ts\eta,\zeta\rangle=\langle \eta,F_{ji}\ts\zeta\rangle,
\quad \eta,\zeta\in V(\lambda),
\non
\end{equation}
cf. \eqref{aform}.
One easily shows that its restriction to
the subspace $V(\lambda)^+_{\mu}$ is non-degenerate. Therefore,
our claim will follow from the fact that the submodules
generated by $\xi_{\mu}$ and $\xi'_{\mu}$ are orthogonal
to each other with respect to $\langle\ts,\rangle$.
Given an operator $A$ in $V(\lambda)^+$ we denote by
$A^*$ its adjoint operator with respect to the form:
\beq%\label{adj}
\langle A\ts\eta,\zeta\rangle=\langle \eta,A^*\ts\zeta\rangle.
\non
\end{equation}
Since the
extremal projection $p$ is stable with respect to the anti-involution
$F_{ij}\mapsto F_{ji}$ \cite{z:it}
we derive that $(pF_{ia})^*=pF_{ai}$ for $a=-n,n$ and 
$i=-n+1,\dots,n-1$. Therefore, $z_{ia}^*=z_{ai}\cdot c$
where $c$ is an element of $\R(\h_{n-1})$ which can be found from
\eqref{zia} and \eqref{zai}. This also implies that
$(z_{ai}z_{ib})^*=z_{bi}z_{ia}$ and hence
\beq\label{Zadj} 
Z_{ab}(u)^*=Z_{ba}(u),
\end{equation}
see \eqref{Zab}. By Proposition~\ref{prop:map} and
the Poincar\'e--Birkhoff--Witt theorem
for the twisted Yangians \cite[Remark~3.14]{mno:yc}, every element
of the $\Y^+(2)$-submodules generated by $\xi_{\mu}$ and  $\xi'_{\mu}$
can be written as a linear combination of vectors of the
following form, respectively:
\beq%\label{vect}
Z_{n,-n}(u_1)\cdots Z_{n,-n}(u_k)\ts\xi_{\mu}
\quad\text{or}\quad
Z_{n,-n}(v_1)\cdots Z_{n,-n}(v_l)\ts\xi'_{\mu},
\non
\end{equation}
where the $u_i$ and $v_i$ are complex parameters.
Therefore, by \eqref{ZZxi} and \eqref{Zadj} to prove that
the submodules are orthogonal it now suffices to show that
$\langle \xi_{\mu},\xi'_{\mu}\rangle=0$. But this follows
from \eqref{zon}.
\endproof

\noindent
{\it Remark.} Using Weyl's formula for the dimension of $V(\lambda)$
one can slightly modify the proof of Theorem~\ref{thm:isom}
so that the branching rule \eqref{bra} would not be used but
follow from the theorem; cf. \cite{m:br, m:wb}.
\endproof

It follows from \eqref{sym} that
the series $s_{n,-n}(u)$ is even in $u$, and so 
is the polynomial $Z_{n,-n}(u)$; see Proposition~\ref{prop:map}.
On the other hand, \eqref{Zab} implies that
$Z_{n,-n}(-g_i)=z_{ni}z_{i,-n}$ for $i=1,\dots,n-1$.
Moreover, $Z_{n,-n}(-g_n)=z_{n,-n}$ which follows from 
\eqref{zn-n}.
Since $Z_{n,-n}(u)$ is a polynomial in $u^2$ of degree $n-1$,
by the Lagrange interpolation formula
$Z_{n,-n}(u)$ can also be given by
\begin{equation}\label{Zn-n}
Z_{n,-n}(u)=\sum_{i=1}^{n}z_{ni}z_{i,-n}\prod_{j=1,\ts j\ne i}^{n}
\frac{u^2-g_j^2}{g_i^2-g_j^2}.
\end{equation}

\noindent
{\it Remark.} To make the above evaluation $Z_{n,-n}(-g_i)$
well-defined we agree to consider the series $Z_{ab}(u)$
with $a,b\in\{-n,n\}$
as elements of the {\it right\/} module over the field
of rational functions in $g_1,\dots,g_n,u$ generated by
monomials in the $z_{ia}$. \endproof

Given $\nu$ such that the conditions \eqref{ineq}
are satisfied, set
$$
\gamma_i=\nu_i-i+1,\qquad l_i=\lambda_{i}-i+1,\qquad i\geq 1
$$
and introduce the vectors
\begin{equation}%\label{xinumu}
\xi_{\nu\mu}=\begin{cases}
\displaystyle{%
\prod_{i=1}^{n}Z_{n,-n}(\gamma_i-1)\cdots 
Z_{n,-n}(\beta_i+1)Z_{n,-n}(\beta_i)\ts\xi_{\mu}}
\qquad&\text{if $\sigma=0$},\\
\displaystyle{%
\prod_{i=1}^{n}Z_{n,-n}(\gamma_i-1)\cdots 
Z_{n,-n}(\beta_i+1)Z_{n,-n}(\beta_i)\ts\xi'_{\mu}}
\qquad&\text{if $\sigma= 1$}.
\end{cases}
\non
\end{equation}
Using \eqref{Zn-n} and \eqref{relcom}
we get an equivalent expression; 
cf. \cite[Section~6]{m:br}:
\begin{equation}%\label{xinumuz}
\xi_{\nu\mu}=
z_{n0}^{\sigma}\ts\prod_{i=1}^{n-1}z_{ni}^{\nu_{i}-\mu_i}
z_{i,-n}^{\nu_{i}-\lambda_i}\cdot
\prod_{k=l_{n}}^{\gamma_{n}-1}Z_{n,-n}(k)\ts\xi.
\non        
\end{equation}

\bpr\label{prop:basis}
The vectors $\xi_{\nu\mu}$ with $\nu$ satisfying \eqref{ineq}
form a basis of $V(\lambda)^+_{\mu}$.
\epr

\Proof Due to Theorem~\ref{thm:isom}
it suffices to show that 
the vectors $\xi_{\nu\mu}$ with $\sigma=0$
form a basis of the subspace $U$ of
$V(\lambda)^+_{\mu}$ while those with $\sigma= 1$ form
a basis in $U'$.
Let us write each of the tensor
products in \eqref{uu}--\eqref{uu2} in the form
\beq\label{genten}
L(\alpha_1,\beta_1)\ot L(\alpha_2,\beta_2)\ot\cdots\ot
L(\alpha_n,\beta_n)\ot W(-\alpha_0).
\end{equation}
Regarding this as a $\Y(2)$-module defined by \eqref{isomL}
we can construct a Gelfand--Tsetlin-type basis in this module as follows.
Set
\beq%\label{GTY}
\wt{\zeta}_{\nu\mu}=
\prod_{i=1}^{n}T_{n,-n}(-\gamma_i+1)\cdots 
T_{n,-n}(-\beta_i-1)T_{n,-n}(-\beta_i)\ts 
z^{\sigma}_{n0}\ts\xi_{\mu}.
\non
\end{equation}
The vectors $\wt{\zeta}_{\nu\mu}$ with $\nu$ satisfying \eqref{ineq}
form a basis in the $\Y(2)$-module \eqref{genten};
see \cite{t:im, m:gt, nt:ry}. 
Furthermore, we have
\begin{equation}\label{TTT}
\begin{split}
T_{nn}(u)\ts\wt{\zeta}_{\nu\mu}&=(u+\gamma_1)\cdots (u+\gamma_{n})
\ts\wt{\zeta}_{\nu\mu}, \\
T_{n,-n}(-\gamma_i)\ts \wt{\zeta}_{\nu\mu}&=
\ts
\wt{\zeta}_{\nu+\delta_i,\mu}.
\end{split}
\end{equation}
We have the following equality of operators in the space
\eqref{genten}:
\begin{equation}%\label{ZT1}
Z_{n,-n}(u)=
\frac{(u-\alpha_0)T_{n,-n}(-u)T_{nn}(u)+
(u+\alpha_0)T_{n,-n}(u)T_{nn}(-u)}{(-1)^{n+1}\ts u}
\non
\end{equation}
which is easily derived from
\eqref{rel}, \eqref{cops} and \eqref{map}. Therefore, by \eqref{TTT}
\begin{equation}
Z_{n,-n}(\gamma_i)\ts\wt{\zeta}_{\nu\mu}=
-2(\alpha_0-\gamma_i)\ts\prod_{a=1,\ts a\ne i}^{n}(-\gamma_a-\gamma_i)
\ts \wt{\zeta}_{\nu+\delta_i,\mu}.
\non
\end{equation}
This shows that for each $\nu$
the vectors $\xi_{\nu\mu}$ and $\wt{\zeta}_{\nu\mu}$ coincide
up to a nonzero factor. \endproof

We shall use the following
normalized basis vectors
\begin{equation}%%\label{zetanumu}
\zeta_{\nu\mu}=\prod_{1\leq i<j\leq n}
(-\gamma_i-\gamma_j)!\  \xi_{\nu\mu}.
\non
\end{equation}
The following formulas for
the action of
the generators of the Yangian $\Y(2)$
in the basis $\{\zeta_{\nu\mu}\}$
follow from the above proof: for $i=1,\dots,n$
\begin{align}
T_{nn}(u)\ts\zeta_{\nu\mu}&=(u+\gamma_1)\cdots (u+\gamma_{n})
\ts\zeta_{\nu\mu}, \non\\
T_{n,-n}(-\gamma_i)\ts \zeta_{\nu\mu}&=
\frac{1}{2(\gamma_i-\alpha_0)}\ts
\zeta_{\nu+\delta_i,\mu},\label{acT1}
\\
T_{-n,n}(-\gamma_i)\ts \zeta_{\nu\mu}&=2\ts
\prod_{k=0}^{n}(\alpha_k-\gamma_i+1)\prod_{k=1}^{n}(\beta_k-\gamma_i)
\ts\zeta_{\nu-\delta_i,\mu};\non
\end{align}
cf. \cite{m:br} and \cite{m:wb}.

Given a pattern $\Lambda$ (see Section~\ref{subsec:mt})
introduce the vector
\begin{equation}%%\label{basisxi}
\xi^{}_{\Lambda}=
\prod_{k=1,\dots,n}^{\rightarrow}
\left(z_{k0}^{\sigma^{}_{k}}\cdot\prod_{i=1}^{k-1}
z_{ki}^{\lambda'_{ki}-\lambda^{}_{k-1,i}}
z_{i,-k}^{\lambda'_{ki}-\lambda^{}_{ki}}
\cdot\prod_{q=l^{}_{kk}}^{l'_{kk}-1}
Z_{k,-k}(q+\frac12)\right)\xi
\non
\end{equation}
and set
\begin{equation}%\label{basiszeta}
\zeta^{}_{\Lambda}=N^{}_{\Lambda}\ts \xi^{}_{\Lambda},\qquad
N^{}_{\Lambda}=\prod_{k=2}^{n}\prod_{1\leq i<j\leq k}(-l'_{ki}-l'_{kj}-1)!
\non
\end{equation}
The following proposition is implied by	the branching rule \eqref{bra}
and Proposition~\ref{prop:basis}.

\bpr\label{prop:totbasis}
The vectors $\zeta^{}_{\Lambda}$ parametrized by the 
%Gelfand--Tsetlin 
patterns $\Lambda$ form a basis of 
the representation $V(\lambda)$. \endproof
\epr

\subsection{Matrix element formulas}\label{subsec:me}

Introduce the
following elements of $\U(\g_n)$:
\begin{equation}%\label{Phi}
\Phi_{-k,k}=\sum_{i=1}^{k-1}F_{-k,i}F_{ik}-\frac12\ts F_{0k}^2,
\qquad k=1,\dots,n.
\non
\end{equation}
We shall find the action of $\Phi_{-k,k}$ in the basis $\{\zeta^{}_{\Lambda}\}$,
which will be used later on.
Since $\Phi_{-k,k}$ commutes with the subalgebra $\g_{k-1}$ it suffices
to consider the case $k=n$. The image of $2\Phi_{-n,n}$ under
the natural homomorphism $\pi:\CC_n\to\Z(\g_n,\g_{n-1})$
coincides with the coefficient at $u^{2n-2}$ of the polynomial
$Z_{-n,n}(u)$; see the proof of Proposition~\ref{prop:map}.
The following equality of operators in \eqref{genten} is obtained from
\eqref{rel}, \eqref{cops} and \eqref{map}:
\begin{equation}%\label{ZT2}
Z_{-n,n}(u)=
\frac{(u-\alpha_0)T_{-n,-n}(-u)T_{-n,n}(u)+
(u+\alpha_0)T_{-n,-n}(u)T_{-n,n}(-u)}{(-1)^{n+1}\ts u}.
\non
\end{equation}
Therefore,
\begin{equation}\label{Phit}
\Phi_{-n,n}=-t^{(2)}_{-n,n}+t^{(1)}_{-n,n}t^{(1)}_{-n,-n}+
(1+\alpha_0)\ts t^{(1)}_{-n,n}.
\end{equation}
The image of $s_{nn}^{(1)}$ under the homomorphism \eqref{homom}
is $F_{nn}$. On the other hand, by \eqref{cops} we have
\begin{equation}%\label{snn1}
s_{nn}^{(1)}=t_{nn}^{(1)}-t_{-n,-n}^{(1)}-\alpha_0-1/2,
\non
\end{equation}
as operators in the space \eqref{genten}. 
Therefore, \eqref{Phit} can be written as
\begin{equation}%\label{Phit2}
\Phi_{-n,n}=-t^{(2)}_{-n,n}+t^{(1)}_{-n,n}t^{(1)}_{nn}-
(F_{nn}+3/2)\ts t^{(1)}_{-n,n}.
\non
\end{equation}
Finally, relations \eqref{acT1} imply that
\begin{equation}\label{Phiact}
\Phi_{-n,n}\zeta_{\nu\mu}=\sum_{i=1}^{n}\theta_i\ts(F_{nn}-\gamma_i+3/2)
\ts \zeta_{\nu-\delta_i,\mu},
\end{equation}
where
\begin{equation}%\label{theta}
\theta_i=-2\ts\prod_{k=0}^{n}(\alpha_k-\gamma_i+1)
\prod_{k=1}^{n}(\beta_k-\gamma_i)
\prod_{j=1,\ts j\ne i}^{n}(\gamma_j-\gamma_i)^{-1}.
\non
\end{equation}
Using the notation in \eqref{genten} we can also write this as
\beq%\label{theta2}
\theta_i=(2\gamma_i-1)(1-\sigma-\gamma_i)
\ts\prod_{k=1}^{n}(l_k-\gamma_i)
\prod_{k=1}^{n-1}(m_k-\gamma_i)
\prod_{j=1,\ts j\ne i}^{n}(\gamma_j-\gamma_i)^{-1}.
\non
\end{equation}
The action of $F_{nn}$ in $V(\lambda)^+_{\mu}$ is immediately found so that
\begin{equation}%\label{Fnn}
F_{nn}\ts\xi_{\nu\mu}=\left(\sigma+2\ts\sum_{i=1}^{n}\nu_i-
\sum_{i=1}^n\lambda_i-\sum_{i=1}^{n-1}\mu_i\right)\xi_{\nu\mu}.
\non
\end{equation}

The operator
$F_{n-1,-n}$ preserves
the subspace of $\g_{n-2}$-highest vectors in $V(\lambda)$.
Therefore it suffices to calculate its action
on the basis vectors of the form
\begin{equation}%\label{numunu}
\xi_{\nu\mu\nu'}=X_{\mu\nu'}\ts\xi_{\nu\mu},
\non
\end{equation}
where 
$X_{\mu\nu'}$ denotes the operator
\begin{equation}%\label{X}
X_{\mu\nu'}=
z_{n-1,0}^{\sigma'}\prod_{i=1}^{n-2}
z_{n-1,i}^{\nu'_{i}-\mu'_i}\ts z_{i,-n+1}^{\nu'_{i}-\mu_i}\cdot
\prod_{a=m_{n-1}}^{\gamma'_{n-1}-1}
Z_{n-1,-n+1}(a).
\non
\end{equation}
Here we assume that
the conditions \eqref{ineq} are satisfied with
$\lambda$, $\nu$, $\mu$ respectively replaced
by $\mu$, $\nu'$, $\mu'$; we have used the notation
$\gamma'_i=\nu'_i-i+1$.
The operator $F_{n-1,-n}$ is permutable with the 
elements $z_{n-1,i}$, $z_{i,-n+1}$ and
$Z_{n-1,-n+1}(u)$ which follows from their explicit formulas. 
Hence, we can write
\begin{equation}%\label{FX}
F_{n-1,-n}\ts \xi_{\nu\mu\nu'}=X_{\mu\nu'}\ts F_{n-1,-n}\ts \xi_{\nu\mu}.
\non
\end{equation}
Let us apply \eqref{Fn-1a} with $a=-n$. We have
\beq%\label{fiact}
f_i\ts \xi_{\nu\mu}=(m_i-1/2)\ts\xi_{\nu\mu},\qquad 
f_{-i}\ts \xi_{\nu\mu}=(-m_i+1/2)\ts\xi_{\nu\mu}
\non
\end{equation}
for $i\geq 1$. Recall also that $f_0=-1/2$.
We now need to express
\begin{equation}%\label{Xgen}
X_{\mu\nu'}\ts z_{n-1,i}z_{i,-n}\ts \xi_{\nu\mu},\qquad
i=-n+1,\dots,n-1
\non
\end{equation}
as a linear combination of the vectors $\xi_{\nu\mu\nu'}$.
Suppose first that $i\geq 1$. Assuming that $\nu_i-\mu_i\geq 1$
we obtain from \eqref{relcom}
\beq\label{ipol}
z_{i,-n}\ts \xi_{\nu\mu}=z_{i,-n}z_{ni}\ts \xi_{\nu,\mu+\delta_i}.
\end{equation}
By \eqref{Zab} this equals
\beq%\label{ziacting}
z_{i,-n}z_{ni}\ts \xi_{\nu,\mu+\delta_i}
=z_{n,-i}z_{-i,-n}\ts \xi_{\nu,\mu+\delta_i}
=Z_{n,-n}(-g_{-i})\ts \xi_{\nu,\mu+\delta_i}.
\non
\end{equation}
We have $-g_{-i}\ts \xi_{\nu,\mu+\delta_i}=m_i\ts \xi_{\nu,\mu+\delta_i}$.
Since 
%\beq%\label{gammap}
$Z_{n,-n}(\gamma_p)\ts \xi_{\nu,\mu+\delta_i}
=\xi_{\nu+\delta_p,\mu+\delta_i}$
%\non
%\end{equation}
for each $p=1,\dots,n$, we obtain from the 
Lagrange interpolation formula that \eqref{ipol} takes the form
\beq%\label{Zmi}
Z_{n,-n}(m_i)\ts \xi_{\nu,\mu+\delta_i}
=\sum_{p=1}^{n}\prod_{a=1,\ts a\ne p}^n
\frac{m_i^2-\gamma_a^2}{\gamma_p^2-\gamma_a^2}
\ts\xi_{\nu+\delta_p,\mu+\delta_i}.
\non
\end{equation}
Furthermore, for $i\geq 1$
\beq%\label{ineg}
z_{-i,-n}\ts \xi_{\nu\mu}=(-1)^{n-i}\ts z_{ni}\ts \xi_{\nu\mu}
=(-1)^{n-i}\ts \xi_{\nu,\mu-\delta_i}.
\non
\end{equation}
Consider now the vector $z^{}_{0,-n}\ts\xi_{\nu\mu}$.
If $\sigma=0$ then it equals
\beq%\label{izero}
z^{}_{0,-n}\ts\xi_{\nu\mu}=(-1)^n\ts z^{}_{n0}\ts\xi_{\nu\mu}
=(-1)^n\ts \xi_{\overline{\nu}\mu},
\non
\end{equation}
where 
%\beq\label{nubar}
$\overline{\nu}=(\sigma+1,\nu_1,\dots,\nu_n)$
%\end{equation}
(addition modulo $2$). If $\sigma=1$ then
\beq%\label{izero2}
z^{}_{0,-n}\ts\xi_{\nu\mu}=z^{}_{n0}z^{}_{0,-n}
\ts\xi_{\overline{\nu}\mu},
\non
\end{equation}
which coincides with $Z_{n,-n}(-g_{0})\ts\xi_{\overline{\nu}\mu}$,
where $g_{0}=0$. Using again the Lagrange interpolation formula
we find that this equals
\beq%\label{Zmzero}
Z_{n,-n}(m_0)\ts \xi_{\overline{\nu}\mu}
=\sum_{p=1}^{n}\prod_{a=1,\ts a\ne p}^n
\frac{m_0^2-\gamma_a^2}{\gamma_p^2-\gamma_a^2}
\ts\xi_{\overline{\nu}+\delta_p,\mu}
\non
\end{equation}
with $m_0=0$. The operator $X_{\mu\nu'}\ts z_{n-1,i}$
is transformed in exactly the same manner; cf. \cite{m:br}.
Combining the results we obtain
\begin{equation}%\label{facti}
F_{n-1,-n}\ts \xi_{\nu\mu\nu'}=A_0\ts\xi(0)+
\sum_{i=1}^{n-1} A_i
\left(
\frac{1}{m_i}\ts \xi^+(i)-\frac{1}{m_i-1}\ts \xi^-(i)\right),
\non
\end{equation}
where
\begin{equation}%\label{ai}
A_i=\prod_{a=1,\ts a\ne i}^{n-1}\frac{1}{m_i-m_a}
\prod_{a=1}^{n-1}\frac{1}{m_i+m_a-1}.
\non
\end{equation}
Furthermore,
\begin{align}%\label{xi+-}
\xi^-(i)&=\xi_{\nu,\mu-\delta_i,\nu'},\non\\
\xi^+(i)&=\sum_{p=1}^{n}\sum_{q=1}^{n-1}
\prod_{a=1,\ts a\ne p}^n\frac{m_i^2-\gamma_a^2}
{\gamma_p^2-\gamma_a^2}\prod_{a=1,\ts a\ne q}^{n-1}\frac{m_i^2-{\gamma'_a}^2}
{{\gamma'_q}^2-{\gamma'_a}^2}\ts 
\xi_{\nu+\delta_p,\mu+\delta_i,\nu'+\delta_q},
\non
\end{align}
and
\begin{alignat}{2}%\label{xi0}
\xi(0)&=(-1)^n\xi_{\overline{\nu}\mu\overline{\nu}'}\quad&&\text{if}\quad
\sigma=\sigma'=0,\non\\
{}&=\sum_{p=1}^{n}
\prod_{a=1,\ts a\ne p}^n\frac{m_0^2-\gamma_a^2}
{\gamma_p^2-\gamma_a^2}\ts 
\xi_{\overline{\nu}+\delta_p,\mu,\overline{\nu}'}\quad&&\text{if}\quad
\sigma=1,\ \sigma'=0,\non\\
{}&=-\sum_{q=1}^{n-1}\prod_{a=1,\ts a\ne q}^{n-1}\frac{m_0^2-{\gamma'_a}^2}
{{\gamma'_q}^2-{\gamma'_a}^2}\ts
\xi_{\overline{\nu},\mu,\overline{\nu}'+\delta_q}\quad&&\text{if}\quad
\sigma=0,\ \sigma'=1,\non\\
{}&=(-1)^{n-1}\sum_{p=1}^{n}\sum_{q=1}^{n-1}
\prod_{a=1,\ts a\ne p}^n\frac{m_0^2-\gamma_a^2}
{\gamma_p^2-\gamma_a^2}\prod_{a=1,\ts a\ne q}^{n-1}\frac{m_0^2-{\gamma'_a}^2}
{{\gamma'_q}^2-{\gamma'_a}^2}\ts
&&\xi_{\overline{\nu}+\delta_p,\mu,\overline{\nu}'+\delta_q}\non\\
&&&\text{if}\quad
\sigma=\sigma'=1, 
\non
\end{alignat}
with
%\beq\label{overnu}
$\overline{\nu}'=(\sigma'+1,\nu'_1,\dots,\nu'_{n-1})$
%\end{equation}
(addition modulo $2$).
We now
compute the action of 
$F_{n-1,n}$. In the formula \eqref{Fn-1a} with $a=n$
replace the operators $z_{in}$ by
with following expression:
for $i=-n+1,\dots,n-1$
\begin{equation}\label{zincom}
z_{in}=[z_{i,-n},\Phi_{-n,n}]\ts \frac{1}{f_i+F_{nn}}
\end{equation}
and then use the formulas for the action of $z_{i,-n}$ and $\Phi_{-n,n}$;
see \eqref{Phiact}. More precisely, we regard \eqref{zincom}
as a relation in the transvector algebra $\Z(\g_n,\g_{n-1})$
which can be proved 
as follows. First, we calculate 
the commutator $[F_{i,-n},\Phi_{-n,n}]$ in $\U(\g_n)$
then consider it modulo the ideal $\J$
and apply the extremal projection $p$
(see Section 2). 

We have
$
\Phi_{-n,n}F_{nn}=(F_{nn}+2)\ts\Phi_{-n,n}
$
and so, \eqref{Fn-1a} and \eqref{zincom} imply that
\begin{equation}\label{Fact}
F_{n-1,n}\ts\xi_{\nu\mu\nu'}
=X_{\mu\nu'}\ts
\left(\Phi_{n-1,-n}(2)\ts\Phi_{-n,n}-
\Phi_{-n,n}\Phi_{n-1,-n}(0)\right)\ts\xi_{\nu\mu},
\end{equation}
where
\begin{equation}\label{Phinn}
\Phi_{n-1,-n}(u)=\sum_{i=-n+1}^{n-1}z_{n-1,i}\ts z_{i,-n}
\prod_{a=-n+1,\ts a\ne i}^{n-1}\frac{1}{f_i-f_a}\cdot\frac{1}{u+f_i+F_{nn}}.
%\non
\end{equation}
The action of $\Phi_{n-1,-n}(u)$
is found exactly as that of $F_{n-1,-n}$.
Formula \eqref{Fact} is valid provided the denominators
in \eqref{Phinn} do not vanish. However, since \eqref{zincom}
holds in the transvector algebra $\Z(\g_n,\g_{n-1})$,
the relation \eqref{Fact} holds for generic
parameters $\nu$, $\mu$ and $\nu'$ which allows one to get
explicit formulas for all matrix elements of $F_{n-1,n}$.

Finally, the proof of Theorem~\ref{thm:main} is completed by
rewriting the above formulas for the action of the generators
in terms of the parameters $\sigma^{}_k$, $l^{}_{ki}$ and $l'_{ki}$
of the patterns $\Lambda$. The parameters $l_i$, $\gamma_i$, $m_i$
are replaced by
\beq%\label{repl}
l_i\mapsto l^{}_{ki}+1/2,\qquad \gamma_i\mapsto l'_{ki}+1/2,\qquad
m_i\mapsto l^{}_{k-1,i}+1/2.
\non
\end{equation}

\medskip
\noindent
{\it Remark\/}. There is another way of
calculating the matrix elements of $F_{k-1,k}$ based on the
formulas for the action of $T_{-n,-n}(u)$ in the basis
$\{\zeta_{\nu\mu}\}$; see \cite{m:wb}.

\end{document}